\documentstyle[12pt]{amsart} \mathsurround=2pt
\newtheorem{thm}{Theorem}[section]
\newtheorem{crl}[thm]{Corollary}
\newtheorem{prp}[thm]{Proposition}

\begin{document}

\title{Novikov operad is not Koszul}
\author{A.S. Dzhumadil'daev}
\address{Institute of Mathematics, Pushkin street 125, Almaty, KAZAKHSTAN }
\email{askar56@@hotmail.com}



\maketitle
\begin{abstract} An algebra with identities $a\circ(b\circ c-c\circ b)=(a\circ b)\circ c-(a\circ c)\circ b$
and $a\circ(b\circ c)=b\circ(a\circ c)$ is called Novikov. We show that Novikov operad is not Koszul.
\end{abstract}

\section{Introduction}

After putting my paper \cite{Dzh2}  on ArXive I receive several questions concerning non-Koszulity of Novikov operad. 
In \cite{Dzh2} I refer to \cite{Dzh3}  to non-Koszulity of Novikov operad. It is mistake. 
 I am gratefull to P. Zusmanovich and V.Dotsenko who have noticed that the 
 paper \cite{Dzh3} does not contain non-Koszulity of Novikov operad.  
By these reasons  I am putting on ArXive proof of this statement. We will save notations of  \cite{Dzh2}  . 

Recall that an algebra $(A,\circ)$ is called {\it right-symmetric} if it satisfies the identity 
$$(a,b,c)=(a,c,b), \qquad \forall a,b,c\in A,$$
where 
$$(a,b,c)=a\circ(b\circ c)-(a\circ b)\circ c$$
is associator. Algebra is {\it left-commutative},  if 
$$a\circ(b\circ c)=b\circ(a\circ c),$$
for any $a,b,c\in A.$ Similarly, an algebra with identity 
$$(a,b,c)=(b,a,c)$$
is {\it left-symmetric.} The identity
$$(a\circ b)\circ c=(a\circ c)\circ b$$
is called {it right-commutative} identity. 

A right-commutative left-symmetric algebra $(A,\circ)$ is called {\it left-Novikov.}
A left-commutative right-symmetric algebra is called {\it right-Novikov.}  About Novikov algebras see 
\cite{Novikov}, \cite{Gelfand}, \cite{Osborn}.

\begin{prp}\label{44} For any  left-(right-)Novikov algebra $(A,\circ)$ a new algebra $(A,\star)$, where 
$a\star b=b\circ a$ is an opposite multiplication,  became right-(left-)Novikov. 
\end{prp}

\begin{crl}\label{444} Let $F_n^{r}$ be a multilinear part of free right-Novikov algebra generated by $n$ elements and similarly let $F_n^{l}$ be a multilinear part of free left-Novikov algebra generated by $n$ elements. Then $\dim F_n^{r}=\dim F_n^{l}.$
\end{crl}

The aim of our paper is to prove the following result. 

\begin{thm} \label{Non-Koz}  
Dual operad to left-(right-)Novikov operad is right-(left)-Novikov.
Novikov operad is not Koszul.
\end{thm}


{\bf Proof.} 
Free magma has 12 multilinear  elements of degree 3.  For a free Novikov algebras only 6 of them form base. Let us construct multilinear base elements of degree $3$ for free Novikov algebras. 
Bases of free Novikov algebras in terms of rooted trees and so-called $r$-elements was constructed in
\cite{Dzh1}. 
We use free Novikov base construction by Young diagrams given in \cite{Dzh2}. 
 
 \medskip
 
{\noindent Young diagrams of degree  $2$:} 

$$\begin{array}{c}\bullet\\  \bullet\\  
\end{array}\qquad
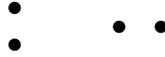
\begin{array}{cc}\bullet&\bullet\\
\end{array}
$$

\medskip

{\noindent Novikov diagrams of degree 3:}

$$
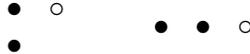
\begin{array}{cc}\bullet&\circ\\
\bullet&\\
\end{array}\qquad
\begin{array}{ccc}\bullet&\bullet&\circ\\
\end{array}
$$

\medskip

{\noindent Novikov tableaux of degree 3 generated by elements $a,b,c$:}

$$
\begin{array}{cc}b&c\\a&\end{array}\qquad
\begin{array}{cc}c&b\\a&\end{array}\qquad
\begin{array}{cc}c&a\\b&\end{array}\qquad
\begin{array}{ccc}a&b&c\\
\end{array}\qquad
\begin{array}{ccc}b&a&c\\
\end{array}\qquad
\begin{array}{ccc}c&a&b\\
\end{array}
$$

\bigskip

{\noindent Multilinear base elements of degree $3$:}
$$a\circ(b\circ c),\quad  a\circ(c\circ b), \quad c\circ(a\circ b), \quad (a\circ b)\circ c, \quad (b\circ a)\circ c, \quad (c\circ a)\circ b.$$ 

Below we give presentation of 6 non-base elements of degree 3 as a linear combination of
above constructed $6$  base elements of degree 3: 
$$b\circ(a\circ c)=a\circ(b\circ c), \quad c\circ(a\circ b)=a\circ(c\circ b), \quad  c\circ(b\circ a)=b\circ (c\circ a),$$ 
$$(a\circ c)\circ b=(a\circ b)\circ c +a\circ(c\circ b) -a\circ(b\circ c),$$
$$(b\circ c)\circ a=-a\circ(b\circ c)+b\circ (c\circ a)+ (b\circ a)\circ c,$$
$$(c\circ b)\circ a=(c\circ a)\circ b-a\circ(c\circ b)+ b\circ (c\circ a).$$

\medskip
Then $$[[a\otimes u,b\otimes v],c\otimes w]$$

$$=((a\circ b)\circ c)\otimes ((u\cdot v)\cdot w)-((b\circ a)\circ
c)\otimes ((v\cdot u)\cdot w)-(c\circ(a\circ b))\otimes( w\cdot(u\cdot
v))$$ $$+(c\circ(b\circ a))\otimes (w\cdot (v\cdot u))$$

$$=((a\circ b)\circ c)\otimes( (u\cdot v)\cdot w) -((b\circ a)\circ c)\otimes ((v\cdot u)\cdot w)$$
$$- (a\circ(c\circ b))\otimes (w\cdot(u\cdot v))+(b\circ (c\circ a))\otimes( w\cdot(v\cdot u))$$ 

$$=((a\circ b)\circ c)\otimes ((u\cdot v)\cdot w)+((b\circ a)\circ c)\otimes (-(v\cdot u)\cdot w)$$
$$+ (a\circ(c\circ b))\otimes(-w\cdot(u\cdot v))+(b\circ (c\circ a))\otimes ( w\cdot (v\cdot u)).$$ 

\medskip

Similarly, $$[[b\otimes v,c\otimes w],a\otimes u]$$

$$=((b\circ c)\circ a)\otimes(( v\cdot w)\cdot u) -((c\circ
b)\circ a)\otimes(( w\cdot v)\cdot u)$$ $$-(a\circ(b\circ c))\otimes
(u\cdot (v\cdot w))+a\circ(c\circ b))\otimes(u\cdot( w\cdot v))$$

$$=(-a\circ(b\circ c)+b\circ (c\circ
a)+(b\circ a)\circ c)\otimes ((v\cdot
w)\cdot u)$$
$$ +(- (c\circ a)\circ b +a\circ(c\circ b) -
b\circ (c\circ a))\otimes(( w\cdot v)\cdot u)$$
$$-(a\circ(b\circ c))\otimes(u\cdot( v\cdot w))+(a\circ(c\circ b))\otimes(u\cdot( w\cdot v))$$

$$=(a\circ(b\circ c))\otimes(-(v\cdot w)\cdot u)+(b\circ (c\circ a))\otimes((v\cdot w)\cdot u)+((b\circ a)\circ c)\otimes ((v\cdot
w)\cdot u)$$
$$ +(- (c\circ a)\circ b +a\circ(c\circ b) -
b\circ (c\circ a))\otimes ((w\cdot v)\cdot u)$$ $$-(a\circ(b\circ c))\otimes(u\cdot( v\cdot w))+(a\circ(c\circ b))\otimes(u\cdot( w\cdot v)).$$ 

\medskip

Further, 
$$[[c\otimes
w,a\otimes u],b\otimes v]$$

$$=((c\circ a)\circ b)\otimes ((w\cdot u)\cdot v)-((a\circ c)\circ
b)\otimes ((u\cdot w)\cdot v) -(b\circ (c\circ a))\otimes (v\cdot
(w\cdot u))$$ $$+(b\circ (a\circ c))\otimes (v\cdot(u\cdot w))$$

$$=((c\circ a)\circ b)\otimes ((w\cdot u)\cdot v)$$
$$+(-(a\circ b)\circ c - a\circ(c\circ b)+a\circ(b\circ c) )\otimes ((u\cdot w)\cdot v)$$
$$+(b\circ (c\circ a) )\otimes (-v\cdot (w\cdot u))$$
$$+( a\circ(b\circ c) )\otimes (v\cdot(u\cdot w))$$

$$=((c\circ a)\circ b)\otimes ((w\cdot u)\cdot v)+((a\circ b)\circ c)\otimes(-(u\cdot w)\cdot v)+(a\circ(c\circ b))\otimes(-(u\cdot w)\cdot v)$$
$$+a\circ(b\circ c) )\otimes ((u\cdot w)\cdot v)+(b\circ (c\circ a) )\otimes (-v\cdot (w\cdot u))
+( a\circ(b\circ c))\otimes (v\cdot(u\cdot w)).$$ 

\medskip

Therefore, $$[[a\otimes u,b\otimes v],c\otimes w] +[[b\otimes
v,c\otimes w],a\otimes u]+[[c\otimes w,a\otimes u],b\otimes v]$$

$$=((a\circ b)\circ c)\otimes \{(u\cdot v)\cdot w  -(u\cdot w)\cdot v\}$$
$$+(a\circ(b\circ c))\otimes\{ -( v\cdot
w)\cdot u -u\cdot( v\cdot w) +v\cdot(u\cdot w)
+(u\cdot w)\cdot v\}$$
$$+ (a\circ(c\circ b))\otimes\{-
w\cdot(u\cdot v)+ (w\cdot v)\cdot u +u\cdot(w\cdot v)
-(u\cdot w)\cdot v\}$$
$$+((b\circ a)\circ c)\otimes\{-(v\cdot u)\cdot
w + (v\cdot w)\cdot u\}$$ 
$$+(b\circ (c\circ a))\otimes $$ $$\{+w\cdot (v\cdot u) + (v\cdot w)\cdot u - (w\cdot v)\cdot u
-v\cdot (w\cdot u)\}$$
$$ +((c\circ a)\circ b) \otimes \{-(w\cdot
v)\cdot u+(w\cdot u)\cdot v\}.$$ 

Thus the Lie-admissibility
condition for $A\otimes U$, where $A$ is free right-Novikov algebra, is equivalent to the following
conditions 
$$(u\cdot v)\cdot w -(u\cdot w)\cdot v =0,$$
$$-( v\cdot w)\cdot u -u\cdot( v\cdot w) +v\cdot(u\cdot
w) +(u\cdot w)\cdot v)=0,$$
$$- w\cdot(u\cdot v)+ (w\cdot v)\cdot u+u\cdot(w\cdot v)
-(u\cdot w)\cdot v=0,$$
$$-(v\cdot u)\cdot w+(v\cdot w)\cdot u=0,$$ 
$$w\cdot (v\cdot u)+(v\cdot w)\cdot u-(w\cdot v)\cdot u
-v\cdot (w\cdot u)=0,$$
$$-(w\cdot v)\cdot u+(w\cdot u)\cdot v=0.$$
Note that all of these identities are consequences of left-symmetric and right-commutative identities. 
So, $(U,\cdot)$ is left-Novikov, if $(A,\circ)$ is right-Novikov. 
Similarly,  $(U,\cdot)$ is right-Novikov, if $(A,\circ)$ is left-Novikov. 
These mean that dual operad to right-(left-)Novikov operad is left-(right-)Novikov operad. 

We have noted that categories of  left-Novikov and right-Novikov algebras are equivalent. 
If we change  left-Novikov multiplication to opposite multiplication we obtain right-Novikov multiplication.  
By Theorem~1.1 \cite{Dzh2}  and by corollary~\ref{444} of proposition~\ref{44}  
dimensions of multilinear parts of free left-(right-)Novikov algebras have dimensions $1,2,6,20,70$ for degrees $1,2,3,4,5.$
 Therefore, beginning parts of Hilbert series of Novikov and dual Novikov operads looks like 
 $$H(t)=H^{!}(t)=-t+t^2-t^3+20t^4/24-70t^5/120+O(t^6).$$
Thus, $$H(H^{!}(t))=t+t^5/6+O(t^6)\ne t.$$
So, by results of 
\cite{Ginz-Kapr}  left-(right-)Novikov operad is not Koszul.
Theorem~\ref{Non-Koz} is proved.

\end{document}